\documentclass{article}
\usepackage{amsmath}
\usepackage{amsfonts}
\usepackage{amssymb}
\usepackage{graphicx}

\newtheorem{theorem}{Theorem}

\newtheorem{corollary}[theorem]{Corollary}

\newtheorem{example}[theorem]{Example}

\newenvironment{proof}[1][Proof]{\textbf{#1.} }{\ \rule{0.5em}{0.5em}}
\input{tcilatex}

\begin{document}

\title{Integer Powers of Certain Complex $(2k+1)-$diagonal Toeplitz Matrices}
\author{Hatice K\"{u}bra Duru\thanks{%
hkduru@selcuk.edu.tr} and Durmu\c{s} Bozkurt\thanks{%
dbozkurt@selcuk.edu.tr} \\
Selcuk University, Science Faculty Department of Mathematics,\\
Turkey}
\maketitle

\begin{abstract}
In this paper, we obtain a general expression for the entries of the $l$th $%
(l\in \mathbb{Z})$ powers of even order $(2k+1)-$diagonal Toeplitz
matrices.\linebreak Additionally, we have the complex factorizations of
Fibonacci\linebreak polynomials.
\end{abstract}

\section{Introduction}

Tridiagonal matrices are used in solving problems in a large variety of
disciplines including mathematics, physics and others. In \cite{Mazilu} it
indicates that using tridiagonal matrices to finding the magnetization
dynamics of outline an approach. Ahmed et al. \cite{Ahmed} enquired some
explicit formulas for the powers of a square matrix. \"{O}tele\c{s} and
Akbulak \cite{Akbulak} gained a general expression for the entries of the
power of complex tridiagonal matrix, in terms of the Chebyshev polynomials
of the first kind. A general expression for the power of complex
persymmetric antitridiagonal Hankel matrices are studied in \cite{Wang}.
Silva \cite{silva} stated integer powers of real anti-tridiagonal matrices
using Chebyshev polynomials. The expressions for powers of complex skew
circulant matrices are found in \cite{Koken}. Duru and Bozkurt \cite{Bozkurt}
are examined powers of tridiagonal matrices. The powers of Toeplitz matrices
are calculated in \cite{Honglin}.

Note that our results in this paper are the more general form of the results
obtained in \cite{Salkuyeh,Duru}. Salkuyeh \cite{Salkuyeh} are calculated
powers of the tridiagonal toeplitz matrices as $k=1$. One can easily see
this in \cite{Duru}, taking $k=2$ in Theorem 1.

This paper is organized as follows: in the next section, we give the
eigenvalues and eigenvectors of $(2k+1)-$ diagonal $n-$ square Toeplitz
matrix. In Section $3$, the $l-$th power of $(2k+1)-$ diagonal Toeplitz
matrix we will get by using the expression $(A_k)^{l}=$ $%
P_k(J_k)^{l}(P_k)^{-1}$ \cite{P.Horn}, where $J_k$ is the Jordan's form of $%
A_k$, $P_k$ is the transforming matrix. In Section $4$, numerical examples
are given. In Section $5$, determinant of $(2k+1)-$ diagonal Toeplitz matrix
obtain complex factorizations for Fibonacci polynomials.

Let $A_k$ be the $(2k+1)-$diagonal $n-$square Toeplitz matrix as following 
\begin{equation}
A_k=\left[ 
\begin{array}{cccccccc}
a & 0 & \cdots & 0 & b & 0 & \cdots & 0 \\ 
0 & a & 0 & \vdots & 0 & b & \cdots & \vdots \\ 
\vdots & 0 & \ddots & 0 & \vdots & \cdots & \cdots & 0 \\ 
0 & \vdots & \cdots & a & 0 & \cdots & 0 & b \\ 
c & 0 & \cdots & 0 & a & 0 & \cdots & 0 \\ 
0 & c & 0 & \vdots & 0 & a & \cdots & \vdots \\ 
\vdots & 0 & \ddots & 0 & \vdots & \cdots & \cdots & 0 \\ 
0 & \vdots & \cdots & c & 0 & \cdots & 0 & a \\ 
&  &  &  &  &  &  & 
\end{array}
\right]  \label{1}
\end{equation}
where $a\in 
\mathbb{R}
$ and $b,c\in 
\mathbb{R}
\setminus \left\{ 0\right\}$.

For example, $A_{5}$ is given as follows: 
\begin{equation*}
A_{5}=\left[ 
\begin{array}{cccccccccc}
a & 0 & 0 & 0 & 0 & b & 0 & 0 & 0 & 0 \\ 
0 & a & 0 & 0 & 0 & 0 & b & 0 & 0 & 0 \\ 
0 & 0 & a & 0 & 0 & 0 & 0 & b & 0 & 0 \\ 
0 & 0 & 0 & a & 0 & 0 & 0 & 0 & b & 0 \\ 
0 & 0 & 0 & 0 & a & 0 & 0 & 0 & 0 & b \\ 
c & 0 & 0 & 0 & 0 & a & 0 & 0 & 0 & 0 \\ 
0 & c & 0 & 0 & 0 & 0 & a & 0 & 0 & 0 \\ 
0 & 0 & c & 0 & 0 & 0 & 0 & a & 0 & 0 \\ 
0 & 0 & 0 & c & 0 & 0 & 0 & 0 & a & 0 \\ 
0 & 0 & 0 & 0 & c & 0 & 0 & 0 & 0 & a \\ 
&  &  &  &  &  &  &  &  & 
\end{array}
\right]
\end{equation*}

\section{Main Results}

\begin{theorem}
Let $A_k$ be $n-$square $(2k+1)-$diagonal Toeplitz matrix as in (\ref{1}).
Then the eigenvalues and eigenvectors of the $A_k$ are

\begin{equation}
\lambda_r=a-2\sqrt{bc}\;cos\left(\frac{kr\pi}{n+k} \right), \; r=\overline{1,%
\frac{n}{k}}  \label{2}
\end{equation}
for $n=2ks\ (s\in\mathbb{N})$ and

\begin{equation}
(P_k)_j=%
\begin{bmatrix}
U_0(\alpha_r) & \overbrace{0 \ldots 0}^{(k-1)} & \mu^{\frac{1}{2}}
U_1(\alpha_r) & \ldots & \overbrace{0 \ldots 0}^{(k-1)} & \mu^{\frac{n-k}{2k}%
} U_{\frac{n-k}{k}}(\alpha_r) & \overbrace{0 \ldots 0}^{(k-1)}%
\end{bmatrix}%
^T,  \label{3}
\end{equation}

$j=1,1+k,1+2k,\ldots,n-k+1; r=\frac{j+k-1}{k}$

\begin{equation}
(P_k)_j=%
\begin{bmatrix}
0 & U_0(\alpha_r) & \overbrace{0 \ldots 0}^{(k-1)} & \mu^{\frac{1}{2}}
U_1(\alpha_r) & \ldots & \overbrace{0 \ldots 0}^{(k-1)} & \mu^{\frac{n-k}{2k}%
} U_{\frac{n-k}{k}}(\alpha_r) & \overbrace{0 \ldots 0}^{(k-2)}%
\end{bmatrix}%
^{T}  \label{4}
\end{equation}

$j=2,2+k,2+2k,\ldots,n-k+2; r=\frac{j+k-2}{k}$

\begin{equation}
(P_k)_j=%
\begin{bmatrix}
0 & 0 & U_0(\alpha_r) & \overbrace{0 \ldots 0}^{(k-1)} & \mu^{\frac{1}{2}}
U_1(\alpha_r) & \ldots & \overbrace{0 \ldots 0}^{(k-1)} & \mu^{\frac{n-k}{2k}%
} U_{\frac{n-k}{k}}(\alpha_r) & \overbrace{0 \ldots 0}^{(k-3)}%
\end{bmatrix}%
^{T}  \label{5}
\end{equation}

$j=3,3+k,3+2k,\ldots,n-k+3; r=\frac{j+k-3}{k}$

$\vdots $ 
\begin{equation}
(P_k)_j=%
\begin{bmatrix}
\overbrace{0 \ldots 0}^{(k-2)} & U_0(\alpha_r) & \overbrace{0 \ldots 0}%
^{(k-1)} & \mu^{\frac{1}{2}} U_1(\alpha_r) & \ldots & \overbrace{0 \ldots 0}%
^{(k-1)} & \mu^{\frac{n-k}{2k}} U_{\frac{n-k}{k}}(\alpha_r) & 0%
\end{bmatrix}%
^{T}  \label{6}
\end{equation}

$j=k-1,2k-1,3k-1,\ldots,n-1; r=\frac{j+1}{k}$

\begin{equation}
(P_k)_j=%
\begin{bmatrix}
\overbrace{0 \ldots 0}^{(k-1)} & U_0(\alpha_r) & \overbrace{0 \ldots 0}%
^{(k-1)} & \mu^{\frac{1}{2}} U_1(\alpha_r) & \ldots & \overbrace{0 \ldots 0}%
^{(k-1)} & \mu^{\frac{n-k}{2k}} U_{\frac{n-k}{k}}(\alpha_r)%
\end{bmatrix}%
^{T}  \label{7}
\end{equation}

$j=k,2k,3k,\ldots,n; r=\frac{j}{k}$\newline
where $\mu=\frac{c}{b}, \alpha_{r}=\frac{\lambda_{r}-a}{2\sqrt{bc}}$ and $%
U_{n}(.)$ is the $n$th degree Chebyshev polynomial of the second kind.
\end{theorem}

\begin{proof}
Let 
\begin{equation}
F_{n}:=\left\vert 
\begin{array}{cccccc}
\lambda -a & b &  &  &  &  \\ 
c & \lambda -a & b &  &  &  \\ 
& c & \lambda -a & b &  &  \\ 
&  & \ddots  & \ddots  & \ddots  &  \\ 
&  &  & c & \lambda -a & b \\ 
&  &  &  & c & \lambda -a%
\end{array}%
\right\vert .  \label{8}
\end{equation}%
For initial conditions $\det \left( F_{0}\right) =1$ and $\det \left(
F_{1}\right) =\lambda -a$, we have 
\begin{equation}
\det \left( F_{n}\right) =\left( \lambda -a\right) \det \left(
F_{n-1}\right) -bc\det \left( F_{n-2}\right) .  \label{9}
\end{equation}%
The solution of difference equation in (\ref{9}) is 
\begin{equation}
\det \left( F_{n}\right) =\left( bc\right) ^{\frac{n}{2}}U_{n}\left( \theta
\right)   \label{10}
\end{equation}%
where $\theta =\frac{\lambda -a}{2\sqrt{bc}}$ and $U_{n}\left( .\right) $ is
the $n$th degree Chebyshev polynomial of the second kind \cite{Mason}: 
\begin{equation*}
U_{n}\left( x\right) =\frac{sin((n+1)arccosx)}{sin(arccosx)}
\end{equation*}%
All the roots of $U_{n}\left( x\right) $ are included in the interval $\left[
-1,1\right] $.\thinspace\ Let 
\begin{equation}
\left\vert \lambda I_{n}-A_{k}\right\vert =\Delta _{A_{k}}\left( \lambda
\right)   \label{11}
\end{equation}%
\newline
and \newline
\begin{equation}
\Delta _{A_{k}}\left( \lambda \right) =\left\{ 
\begin{array}{lll}
\det (F_{\frac{n}{2}})^{2}-\det (F_{\frac{n}{2}-1})^{2},\; & if & \;k=1\; \\ 
\det (F_{t})^{n-(t-1)k}\det (F_{t-1})^{tk-n},\; & if & n\leq tk,\;1<k<\frac{n%
}{2}\;and\;t\in 
\mathbb{Z}
\\ 
\det (F_{2})^{n-k}\det (F_{1})^{2k-n},\; & if & \;\frac{n}{2}\leq k\leq
n-1.\;%
\end{array}%
\right.   \label{12}
\end{equation}%
So, the equality (\ref{12}) is written as 
\begin{equation}
\Delta _{A_{k}}\left( \lambda \right) =\left\{ 
\begin{array}{lll}
(bc)^{\frac{n}{2}}[U_{\frac{n}{2}}^{2}(\theta )-U_{\frac{n}{2}-1}^{2}(\theta
)],\; & if & \;k=1\; \\ 
(bc)^{\frac{n}{2}}U_{t}^{n-(t-1)k}(\theta )U_{t-1}^{tk-n}(\theta ),\; & if & 
n\leq tk,\;1<k<\frac{n}{2}\;and\;t\in 
\mathbb{Z}
\\ 
(bc)^{\frac{n}{2}}U_{2}^{n-k}(\theta )U_{1}^{2k-n}(\theta ),\; & if & \;%
\frac{n}{2}\leq k\leq n-1\;%
\end{array}%
\right.   \label{13}
\end{equation}%
where $\theta =\frac{\lambda -a}{2\sqrt{bc}}$. There are the relations
between $V_{n}\left( .\right) $, $W_{n}\left( .\right) $ and $U_{n}\left(
.\right) $ polynomials as following

\begin{equation}
\begin{array}{ll}
V_{\frac{n}{2}}(\theta) & =U_{\frac{n}{2}}(\theta)-U_{\frac{n}{2}-1} (\theta)
\\ 
W_{\frac{n}{2}}(\theta) & =U_{\frac{n}{2}}(\theta)+U_{\frac{n}{2}-1} (\theta)%
\end{array}
\label{134}
\end{equation}
here $V_{n}\left(.\right)$ and $W_{n}\left(.\right)$ are the $n$th degree
Chebyshev polynomial of the third and fourth kind, respectively \cite{Mason}%
. Substituting (\ref{134}) into (\ref{13}), we possess 
\begin{equation}
\Delta_{A_k}\left(\lambda\right)=\left\{ 
\begin{array}{lll}
(bc)^{\frac{n}{2}}V_{\frac{n}{2}}(\theta)W_{\frac{n}{2}}(\theta),\; & if & 
\;k=1\; \\ 
(bc)^{\frac{n}{2}}U_{t}^{n-(t-1)k} (\theta) U_{t-1}^{tk-n} (\theta),\; & if
& n\leq tk,\;1<k<\frac{n}{2} \; and \; t \in 
\mathbb{Z}
\\ 
(bc)^{\frac{n}{2}}U_{2}^{n-k} (\theta) U_{1}^{2k-n} (\theta),\; & if & \;%
\frac{n}{2}\leq k\leq n-1.\;%
\end{array}%
\right.  \label{135}
\end{equation}
The eigenvalues of $A_k$ obtained as 
\begin{equation*}
\lambda_r=a-2\sqrt{bc}\;cos\left(\frac{kr\pi}{n+k} \right), for \; r=%
\overline{1,\frac{n}{k}}. 
\end{equation*}
The multiplicity of all the eigenvalues $\lambda _{r}$ ($r=1,2,\ldots ,\frac{%
n}{k}$) of the matrix $A_k$ are $k$. Since $rank(\lambda _{r}I_{n}-A_k)=n-k$%
, for each eigenvalue $\lambda _{r}$ correspond $\frac{n}{k}$ Jordan cells $%
J_{r}(\lambda _{r})$\ in the matrix $J_{k}$. That is%
\begin{equation}
J_{k}=diag(\overbrace{\lambda _{1},\ldots,\lambda_{1}}^{n/k} ,\overbrace{%
\lambda _{2},\ldots,\lambda _{2}}^{n/k},\ldots ,\overbrace{\lambda _{\frac{n%
}{2}},\ldots,\lambda _{\frac{n}{2}}}^{n/k}).  \label{14}
\end{equation}
Consider the relations $(P_k)^{-1} A_k P_k=J_{k}$ \cite{P.Horn} we have to
obtain the matrices $P_k$ and $(P_k)^{-1}$ and derive the expression of the
matrix $(A_k)^{l}$ for $l\in \mathbb{N}$. Let us denote $j$-th column of $P_k
$ by $(P_k)_{j}$ ($j=1,\ldots ,n).$ Then 
\begin{equation}
A_k P_k=( (P_k)_{1}\lambda _{1}\,\ldots \, (P_k)_{\frac{n}{k}} \lambda _{1}
\, \ldots \, (P_k)_{\frac{(n-1)n}{k}} \lambda _{\frac{n}{2}} \, \ldots \,
(P_k)_{n}\lambda _{\frac{n}{2}}).  \label{15}
\end{equation}%
From Eq. (\ref{15}), we have the system of linear equations as follows:

\begin{equation}
\left. 
\begin{array}{ccc}
A_{k} (P_k)_{1} & = & (P_k)_{1}\lambda _{1} \\ 
& \vdots &  \\ 
A_{k} (P_k)_{\frac{n}{k}} & = & (P_k)_{\frac{n}{k}}\lambda _{1} \\ 
A_{k} (P_k)_{\frac{n}{k}+1} & = & (P_k)_{\frac{n}{k}+1}\lambda _{2} \\ 
& \vdots &  \\ 
A_{k} (P_k)_{\frac{2n}{k}} & = & (P_k)_{\frac{2n}{k}}\lambda _{2} \\ 
& \vdots &  \\ 
A_{k} (P_k)_{\frac{(n-1)n}{k}} & = & (P_k)_{\frac{(n-1)n}{k}}\lambda _{\frac{%
n}{2}} \\ 
& \vdots &  \\ 
A_{k} (P_k)_{n} & = & (P_k)_{n}\lambda _{\frac{n}{2}}%
\end{array}%
\right.  \label{16}
\end{equation}
Solving the set of the system of linear equations in (\ref{16}) as regards

\begin{equation*}
(P_k)_j=%
\begin{bmatrix}
U_0(\alpha_r) & \overbrace{0 \ldots 0}^{(k-1)} & \mu^{\frac{1}{2}}
U_1(\alpha_r) & \ldots & \overbrace{0 \ldots 0}^{(k-1)} & \mu^{\frac{n-k}{2k}%
} U_{\frac{n-k}{k}}(\alpha_r) & \overbrace{0 \ldots 0}^{(k-1)}%
\end{bmatrix}%
^{T}
\end{equation*}

$j=1,1+k,1+2k,\ldots,n-k+1; r=\frac{j+k-1}{k}$

\begin{equation*}
(P_n^k)_j=%
\begin{bmatrix}
0 & U_0(\alpha_r) & \overbrace{0 \ldots 0}^{(k-1)} & \mu^{\frac{1}{2}}
U_1(\alpha_r) & \ldots & \overbrace{0 \ldots 0}^{(k-1)} & \mu^{\frac{n-k}{2k}%
} U_{\frac{n-k}{k}}(\alpha_r) & \overbrace{0 \ldots 0}^{(k-2)}%
\end{bmatrix}%
^{T}
\end{equation*}

$j=2,2+k,2+2k,\ldots,n-k+2; r=\frac{j+k-2}{k}$

\begin{equation*}
(P_n^k)_j=%
\begin{bmatrix}
0 & 0 & U_0(\alpha_r) & \overbrace{0 \ldots 0}^{(k-1)} & \mu^{\frac{1}{2}}
U_1(\alpha_r) & \ldots & \overbrace{0 \ldots 0}^{(k-1)} & \mu^{\frac{n-k}{2k}%
} U_{\frac{n-k}{k}}(\alpha_r) & \overbrace{0 \ldots 0}^{(k-3)}%
\end{bmatrix}%
^{T}
\end{equation*}

$j=3,3+k,3+2k,\ldots,n-k+3; r=\frac{j+k-3}{k}$

$\vdots $

\begin{equation*}
(P_n^k)_j=%
\begin{bmatrix}
\overbrace{0 \ldots 0}^{(k-2)} & U_0(\alpha_r) & \overbrace{0 \ldots 0}%
^{(k-1)} & \mu^{\frac{1}{2}} U_1(\alpha_r) & \ldots & \overbrace{0 \ldots 0}%
^{(k-1)} & \mu^{\frac{n-k}{2k}} U_{\frac{n-k}{k}}(\alpha_r) & 0%
\end{bmatrix}%
^{T}
\end{equation*}

$j=k-1,2k-1,3k-1,\ldots,n-1; r=\frac{j+1}{k}$

\begin{equation*}
(P_n^k)_j=%
\begin{bmatrix}
\overbrace{0 \ldots 0}^{(k-1)} & U_0(\alpha_r) & \overbrace{0 \ldots 0}%
^{(k-1)} & \mu^{\frac{1}{2}} U_1(\alpha_r) & \ldots & \overbrace{0 \ldots 0}%
^{(k-1)} & \mu^{\frac{n-k}{2k}} U_{\frac{n-k}{k}}(\alpha_r)%
\end{bmatrix}%
^{T}
\end{equation*}

$j=k,2k,3k,\ldots,n; r=\frac{j}{k}$\newline
where $\mu=\frac{c}{b}, \alpha_{r}=\frac{\lambda_{r}-a}{2\sqrt{bc}}$ and $%
U_{n}(.)$ is the $n$th degree Chebyshev polynomial of the second kind.
\end{proof}

\section{The integer powers of the matrix $A_{k}$}

Considering (\ref{3}), (\ref{4}), (\ref{5}), (\ref{6}) and (\ref{7}), we
write down the matrix $P_k$ 
\begin{equation*}
P_k=\left[ 
\begin{array}{ccccc}
U_{0}(\alpha_1) & 0 & 0 & \cdots & U_{0}(\alpha_2) \\ 
0 & U_{0}(\alpha_1) & 0 & \cdots & 0 \\ 
\vdots & \vdots & \vdots & \vdots & \vdots \\ 
0 & 0 & 0 & \cdots & U_{0}(\alpha_1) \\ 
\mu^{\frac{1}{2}}U_{1}(\alpha_1) & 0 & 0 & \cdots & \mu^{\frac{1}{2}%
}U_{1}(\alpha_2) \\ 
0 & \mu^{\frac{1}{2}}U_{1}(\alpha_1) & 0 & \cdots & 0 \\ 
\vdots & \vdots & \vdots & \vdots & \vdots \\ 
0 & 0 & 0 & \cdots & \mu^{\frac{1}{2}}U_{1}(\alpha_1) \\ 
\vdots & \vdots & \vdots & \vdots & \vdots \\ 
\mu^{\frac{n-k}{2k}}U_{\frac{n-k}{k}}(\alpha_1) & 0 & 0 & \cdots & \mu^{%
\frac{n-k}{2k}}U_{\frac{n-k}{k}}(\alpha_2) \\ 
0 & \mu^{\frac{n-k}{2k}}U_{\frac{n-k}{k}}(\alpha_1) & 0 & \cdots & 0 \\ 
\vdots & \vdots & \vdots & \vdots & \vdots \\ 
0 & 0 & 0 & \cdots & \mu^{\frac{n-k}{2k}}U_{\frac{n-k}{k}}(\alpha_1) \\ 
&  &  &  & 
\end{array}
\right.
\end{equation*}
\begin{equation}
\left. 
\begin{array}{cccccc}
\qquad & \cdots & U_{0}(\alpha_{\frac{n}{k}}) & 0 & \cdots & 0 \\ 
\qquad & \cdots & 0 & U_{0}(\alpha_{\frac{n}{k}}) & \cdots & 0 \\ 
\qquad & \ddots & \vdots & \vdots & \vdots & \vdots \\ 
\qquad & \cdots & 0 & 0 & \cdots & U_{0}(\alpha_{\frac{n}{k}}) \\ 
\qquad & \cdots & \mu^{\frac{1}{2}}U_{1}(\alpha_{\frac{n}{k}}) & 0 & \cdots
& 0 \\ 
\qquad & \cdots & 0 & \mu^{\frac{1}{2}}U_{1}(\alpha_{\frac{n}{k}}) & \cdots
& 0 \\ 
\qquad & \ddots & \vdots & \vdots & \vdots & \vdots \\ 
\qquad & \cdots & 0 & 0 & \cdots & \mu^{\frac{1}{2}}U_{1}(\alpha_{\frac{n}{k}%
}) \\ 
\qquad & \ddots & \vdots & \vdots & \vdots & \vdots \\ 
\qquad & \cdots & \mu^{\frac{n-k}{2k}}U_{\frac{n-k}{k}}(\alpha_{\frac{n}{k}})
& 0 & \cdots & 0 \\ 
\qquad & \cdots & 0 & \mu^{\frac{n-k}{2k}}U_{\frac{n-k}{k}}(\alpha_{\frac{n}{%
k}}) & \cdots & 0 \\ 
\qquad & \ddots & \vdots & \vdots & \vdots & \vdots \\ 
\qquad & \cdots & 0 & 0 & \cdots & \mu^{\frac{n-k}{2k}}U_{\frac{n-k}{k}%
}(\alpha_{\frac{n}{k}}) \\ 
&  &  &  &  & 
\end{array}
\right].  \label{17}
\end{equation}

Now let us find the inverse matrix $(P_k)^{-1}$ of the matrix $P_k$. If we
denote $i$th row of the inverse matrix $(P_k)^{-1}$ by $(P_k)^{-1}_{i}$,
then we have

\begin{equation}
(P_k)_i=%
\begin{bmatrix}
\eta_r U_0(\alpha_r) & \overbrace{0 \ldots 0}^{(k-1)} & \eta_r \mu^{-\frac{1%
}{2}} U_1(\alpha_r) & \ldots & \eta_r \mu^{-\frac{n-k}{2k}} U_{\frac{n-k}{k}%
}(\alpha_r) & \overbrace{0 \ldots 0}^{(k-1)}%
\end{bmatrix}%
^{T},  \label{18}
\end{equation}
$i=1,1+k,1+2k,\ldots,n-k+1; r=\frac{i+k-1}{k}$

\begin{equation}
(P_k)_i=%
\begin{bmatrix}
0 & \eta_r U_0(\alpha_r) & \overbrace{0 \ldots 0}^{(k-1)} & \ldots & \eta_r
\mu^{-\frac{n-k}{2k}} U_{\frac{n-k}{k}}(\alpha_r) & \overbrace{0 \ldots 0}%
^{(k-2)}%
\end{bmatrix}%
^{T},  \label{19}
\end{equation}
$i=2,2+k,2+2k,\ldots,n-k+2; r=\frac{i+k-2}{k}$

\begin{equation}
(P_k)_i=%
\begin{bmatrix}
0 & 0 & \eta_r U_0(\alpha_r) & \overbrace{0 \ldots 0}^{(k-1)} & \ldots & 
\overbrace{0 \ldots 0}^{(k-1)} & \eta_r \mu^{-\frac{n-k}{2k}} U_{\frac{n-k}{k%
}}(\alpha_r) & \overbrace{0 \ldots 0}^{(k-3)}%
\end{bmatrix}%
^{T},  \label{20}
\end{equation}
$i=3,3+k,3+2k,\ldots,n-k+3; r=\frac{i+k-3}{k}$

$\vdots $

\begin{equation}
(P_k)_i=%
\begin{bmatrix}
\overbrace{0 \ldots 0}^{(k-2)} & \eta_r U_0(\alpha_r) & \overbrace{0 \ldots 0%
}^{(k-1)} & \ldots & \overbrace{0 \ldots 0}^{(k-1)} & \eta_r \mu^{-\frac{n-k%
}{2k}} U_{\frac{n-k}{k}}(\alpha_r) & 0%
\end{bmatrix}%
^{T},  \label{21}
\end{equation}
$i=k-1,2k-1,3k-1,\ldots,n-1; r=\frac{i+1}{k}$

\begin{equation}
(P_k)_i=%
\begin{bmatrix}
\overbrace{0 \ldots 0}^{(k-1)} & \eta_r U_0(\alpha_r) & \overbrace{0 \ldots 0%
}^{(k-1)} & \ldots & \overbrace{0 \ldots 0}^{(k-1)} & \eta_r \mu^{-\frac{n-k%
}{2k}} U_{\frac{n-k}{k}}(\alpha_r)%
\end{bmatrix}%
^{T},  \label{22}
\end{equation}
$i=k,2k,3k,\ldots,n; r=\frac{i}{k}$\newline
where $\mu=\frac{c}{b}, \eta_r=\frac{k(4-4\alpha_r^2)}{2(n+k)}$ and $%
\alpha_{r}=\frac{\lambda_{r}-a}{2\sqrt{bc}}$ for $r=\overline{1,\frac{n}{k}}$%
. Thus, we obtain

\begin{equation*}
(P_k)^{-1}=\left[ 
\begin{array}{cccccc}
\eta_1 U_{0}(\alpha_1) & 0 & 0 & \cdots & \eta_1 \mu^{-\frac{1}{2}}
U_1(\alpha_1) & \cdots \\ 
0 & \eta_1 U_{0}(\alpha_1) & 0 & \cdots & 0 & \cdots \\ 
\vdots & \vdots & \vdots & \vdots & \vdots & \ddots \\ 
0 & 0 & 0 & \cdots & \eta_1 U_{0}(\alpha_1) & \cdots \\ 
\eta_2 U_{0}(\alpha_2) & 0 & 0 & \cdots & \eta_2 \mu^{-\frac{1}{2}}
U_{1}(\alpha_2) & \cdots \\ 
0 & \eta_2 U_{0}(\alpha_2) & 0 & \cdots & 0 & \cdots \\ 
\vdots & \vdots & \vdots & \vdots & \vdots & \ddots \\ 
0 & 0 & 0 & \cdots & \eta_2 U_{0}(\alpha_2) & \cdots \\ 
\vdots & \vdots & \vdots & \vdots & \vdots & \cdots \\ 
\eta_{\frac{n}{k}} U_{0}(\alpha_\frac{n}{k}) & 0 & 0 & \cdots & \eta_\frac{n%
}{k} \mu^{-\frac{1}{2}} U_1(\alpha_\frac{n}{k}) & \cdots \\ 
0 & \eta_{\frac{n}{k}} U_{0}(\alpha_\frac{n}{k}) & 0 & \cdots & 0 & \cdots
\\ 
\vdots & \vdots & \vdots & \vdots & \vdots & \cdots \\ 
0 & 0 & 0 & \cdots & \eta_{\frac{n}{k}} U_{0}(\alpha_\frac{n}{k}) & \cdots
\\ 
&  &  &  &  & 
\end{array}
\right.
\end{equation*}
\begin{equation}
\left. 
\begin{array}{ccccc}
\eta_1 \mu^{-\frac{n-k}{2k}} U_{\frac{n-k}{k}} (\alpha_1) & 0 & \cdots & 0 & 
\\ 
0 & \eta_1 \mu^{-\frac{n-k}{2k}} U_{\frac{n-k}{k}} (\alpha_1) & \cdots & 0 & 
\\ 
\vdots & \vdots & \vdots & \vdots &  \\ 
0 & 0 & \cdots & \eta_1 \mu^{-\frac{n-k}{2k}} U_{\frac{n-k}{k}} (\alpha_1) & 
\\ 
\eta_2 \mu^{-\frac{n-k}{2k}-1} U_{\frac{n-k}{k}-1} (\alpha_2) & 0 & \cdots & 
0 &  \\ 
0 & \eta_2 \mu^{-\frac{n-k}{2k}-1} U_{\frac{n-k}{k}-1} (\alpha_2) & \cdots & 
0 &  \\ 
\vdots & \vdots & \vdots & \vdots &  \\ 
0 & 0 & \cdots & \eta_2 \mu^{-\frac{n-k}{2k}-1} U_{\frac{n-k}{k}-1}
(\alpha_2) &  \\ 
\vdots & \vdots & \vdots & \vdots &  \\ 
\eta_\frac{n}{k} \mu^{-\frac{n-k}{2k}} U_{\frac{n-k}{k}} (\alpha_\frac{n}{k})
& 0 & \cdots & 0 &  \\ 
0 & \eta_\frac{n}{k} \mu^{-\frac{n-k}{2k}} U_{\frac{n-k}{k}} (\alpha_\frac{n%
}{k}) & \cdots & 0 &  \\ 
\vdots & \vdots & \vdots & \vdots &  \\ 
0 & 0 & \cdots & \eta_\frac{n}{k} \mu^{-\frac{n-k}{2k}} U_{\frac{n-k}{k}}
(\alpha_\frac{n}{k}) &  \\ 
&  &  &  & 
\end{array}
\right].  \label{23}
\end{equation}

By combining (\ref{14}), (\ref{17}) and (\ref{23}) and using the equality $%
(A_k)^{l}=$ $P_k(J_k)^{l}(P_k)^{-1}$ \cite{P.Horn}, we compute the $l$th
powers of the matrix $A_k$ as \newline
\begin{equation}
(A_k)^{l}=P_k(J_k)^{l}(P_k)^{-1}=Q\left( l\right) =\left( q_{ij}\left(
l\right) \right) .  \label{24}
\end{equation}%
\newline
So, for $i,j=\overline{1,n}$%
\begin{equation}
q_{ij}\left( l\right) =\left\{ 
\begin{array}{lll}
0, & if & other, \\ 
\sum\limits_{r=1}^{\frac{n}{k}} \lambda _{r}^{l} \eta_{r}\mu ^{\frac{i-j}{2k}%
}U_{\frac{i-\varphi_{ij}}{k}}\left( \alpha _{r}\right) U_{\frac{%
j-\varphi_{ij}}{k}}\left( \alpha _{r}\right) , & if & i=j\; and \; i=j+k,j-k,%
\end{array}%
\right.  \label{25}
\end{equation}

\begin{equation}
\varphi_{ij}=\left\{%
\begin{array}{lllll}
k & ; & \;i,j & \equiv & 0 \mod(k)~, \\ 
m & ; & \;i,j & \equiv & m \mod(k)~, \\ 
&  &  &  & 
\end{array}%
\right.  \label{26}
\end{equation}
where $m=\overline{1,k-1}$.

\begin{corollary}
Let $A_{k}$ be $n-$square $(n=2ks,\;s\in \mathbb{N};\;a,b,c\in 
\mathbb{R}
\setminus \left\{ 0\right\} )$\linebreak $(2k+1)-$diagonal Toeplitz matrix
as in (\ref{1}), from Theorem 1 
\begin{equation}
a\neq 2\sqrt{bc}\;cos\left( \frac{kr\pi }{n+k}\right)   \label{266}
\end{equation}%
$\left( r=\overline{1,\frac{n}{k}}\right) $. In that case, there exists the
inverse and negative integer powers of the matrix $A_{k}$.
\end{corollary}

\section{Numerical Examples}

\begin{example}
Setting $n=8$, for $k=1$ and $l=3$, we have 
\begin{equation*}
\begin{array}{lll}
J_{1} & = & diag(\lambda _{1},\lambda _{2},\lambda _{3},\lambda _{4},\lambda
_{5},\lambda _{6},\lambda _{7},\lambda _{8}) \\ 
& = & diag(a-1.879\sqrt{bc},a-1.532\sqrt{bc},a-\sqrt{bc},a-0.347\sqrt{bc},
\\ 
&  & \qquad a+0.347\sqrt{bc},a+\sqrt{bc},a+1.532\sqrt{bc},a+1.879\sqrt{bc})%
\end{array}%
\end{equation*}%
and \newline
$%
\begin{array}{llllll}
(A_{1})^{3} & = & P_{1}(J_{1})^{3}(P_{1})^{-1} & = & Q(3) & =%
\end{array}%
$%
\begin{equation*}
\begin{array}{ll}
(q_{ij}(3)) & =\left[ 
\begin{array}{cccc}
a^{3}+3abc & 3a^{2}b+2b^{2}c & 3ab^{2} & b^{3} \\ 
3a^{2}c+2bc^{2} & a^{3}+6abc & 3a^{2}b+3b^{2}c & 3ab^{2} \\ 
3ac^{2} & 3a^{2}c+3bc^{2} & a^{3}+6abc & 3a^{2}b+3b^{2}c \\ 
c^{3} & 3ac^{2} & 3a^{2}c+3bc^{2} & a^{3}+6abc \\ 
0 & c^{3} & 3ac^{2} & 3a^{2}c+3bc^{2} \\ 
0 & 0 & c^{3} & 3ac^{2} \\ 
0 & 0 & 0 & c^{3} \\ 
0 & 0 & 0 & 0 \\ 
&  &  & 
\end{array}%
\right. 
\end{array}%
\end{equation*}%
\begin{equation*}
\left. 
\begin{array}{cccc}
0 & 0 & 0 & 0 \\ 
b^{3} & 0 & 0 & 0 \\ 
3ab^{2} & b^{3} & 0 & 0 \\ 
3a^{2}b+3b^{2}c & 3ab^{2} & b^{3} & 0 \\ 
a^{3}+6abc & 3a^{2}b+3b^{2}c & 3ab^{2} & b^{3} \\ 
3a^{2}c+3bc^{2} & a^{3}+6abc & 3a^{2}b+3b^{2}c & 3ab^{2} \\ 
3ac^{2} & 3a^{2}c+3bc^{2} & a^{3}+6abc & 3a^{2}b+2b^{2}c \\ 
c^{3} & 3ac^{2} & 3a^{2}c+2bc^{2} & a^{3}+3abc \\ 
&  &  & 
\end{array}%
\right] .
\end{equation*}

For $k=2,l=4,a=i-1,b=i+4$ and $c=i-2$

\begin{equation*}
\begin{array}{lll}
J_{2} & = & diag(\lambda _{1},\lambda _{1},\lambda _{2},\lambda _{2},\lambda
_{3},\lambda _{3},\lambda _{4},\lambda _{4}) \\ 
& = & diag(-1.536-3.884i,-1.536-3.884i,-1.205-0.865i,-1.205-0.865i, \\ 
&  & \qquad -0.795+2.865i,-0.795+2.865i,-0.464+5.884i,-0.464+5.884i)%
\end{array}%
\end{equation*}

and \newline
$%
\begin{array}{llllll}
(A_{2})^{4} & = & P_{2}(J_{2})^{4}(P_{2})^{-1} & = & Q(l) & =%
\end{array}%
$%
\begin{equation*}
\begin{array}{ll}
(q_{ij}(l)) & =\left[ 
\begin{array}{cccc}
174+36i & 0 & 336-256i & 0 \\ 
0 & 174+36i & 0 & 336-256i \\ 
-48+224i & 0 & 429+36i & 0 \\ 
0 & -48+224i & 0 & 429+36i \\ 
-105+90i & 0 & -60+340i & 0 \\ 
0 & -105+90i & 0 & -60+340i \\ 
-36-52i & 0 & -105+90i & 0 \\ 
0 & -36-52i & 0 & -105+90i%
\end{array}%
\right. 
\end{array}%
\end{equation*}%
\begin{equation*}
\left. 
\begin{array}{cccc}
-357-306i & 0 & -396+20i & 0 \\ 
0 & -357-306i & 0 & -396+20i \\ 
492-404i & 0 & -357-306i & 0 \\ 
0 & 492-404i & 0 & -357-306i \\ 
429+36i & 0 & 336-256i & 0 \\ 
0 & 429+36i & 0 & 336-256i \\ 
-48+224i & 0 & 174+36i & 0 \\ 
0 & -48+224i & 0 & 174+36i%
\end{array}%
\right] .
\end{equation*}

For $k=4$%
\begin{equation*}
\begin{array}{lll}
J_{4} & = & diag(\lambda _{1},\lambda _{1},\lambda _{1},\lambda _{1},\lambda
_{2}\lambda _{2},\lambda _{2},\lambda _{2})%
\end{array}%
\end{equation*}

where $\lambda _{1}=a-\sqrt{bc},\lambda _{2}=a+\sqrt{bc}$ and

$%
\begin{array}{llllll}
(A_{4})^{l} & = & P_{4}(J_{4})^{l}(P_{4})^{-1} & = & Q(l) & =%
\end{array}%
$ \newline
\begin{equation*}
\begin{array}{ll}
(q_{ij}(l)) & =\left[ 
\begin{array}{cccccccc}
q_{11}(l) & q_{12}(l) & q_{13}(l) & q_{14}(l) & q_{15}(l) & q_{16}(l) & 
q_{17}(l) & q_{18}(l) \\ 
q_{21}(l) & q_{22}(l) & q_{23}(l) & q_{24}(l) & q_{25}(l) & q_{26}(l) & 
q_{27}(l) & q_{28}(l) \\ 
q_{31}(l) & q_{32}(l) & q_{33}(l) & q_{34}(l) & q_{35}(l) & q_{36}(l) & 
q_{37}(l) & q_{38}(l) \\ 
q_{41}(l) & q_{42}(l) & q_{43}(l) & q_{44}(l) & q_{45}(l) & q_{46}(l) & 
q_{47}(l) & q_{48}(l) \\ 
q_{51}(l) & q_{52}(l) & q_{53}(l) & q_{54}(l) & q_{55}(l) & q_{56}(l) & 
q_{57}(l) & q_{58}(l) \\ 
q_{61}(l) & q_{62}(l) & q_{63}(l) & q_{64}(l) & q_{65}(l) & q_{66}(l) & 
q_{67}(l) & q_{68}(l) \\ 
q_{71}(l) & q_{72}(l) & q_{73}(l) & q_{74}(l) & q_{75}(l) & q_{76}(l) & 
q_{17}(l) & q_{78}(l) \\ 
q_{81}(l) & q_{82}(l) & q_{83}(l) & q_{84}(l) & q_{85}(l) & q_{86}(l) & 
q_{87}(l) & q_{88}(l)%
\end{array}%
\right] 
\end{array}%
\end{equation*}%
$%
\begin{array}{ll}
q_{11}(l) & =q_{22}(l)=q_{33}(l)=q_{44}(l)=q_{55}(l)=q_{66}(l) \\ 
& =q_{77}(l)=q_{88}(l)=0.5(a+\sqrt{bc})^{l}+0.5(a-\sqrt{bc})^{l} \\ 
q_{15}(l) & =q_{26}(l)=q_{37}(l)=q_{48}(l)=[0.5(a+\sqrt{bc})^{l}-0.5(a-\sqrt{%
bc})^{l}]\mu ^{-\frac{1}{2}} \\ 
q_{51}(l) & =q_{62}(l)=q_{73}(l)=q_{84}(l)=[0.5(a+\sqrt{bc})^{l}-0.5(a-\sqrt{%
bc})^{l}]\mu ^{\frac{1}{2}} \\ 
q_{12}(l) & =q_{13}(l)=q_{14}(l)=q_{16}(l)=q_{17}(l)=q_{18}(l)=0 \\ 
q_{21}(l) & =q_{23}(l)=q_{24}(l)=q_{25}(l)=q_{27}(l)=q_{28}(l)=0 \\ 
q_{31}(l) & =q_{32}(l)=q_{34}(l)=q_{35}(l)=q_{36}(l)=q_{38}(l)=0 \\ 
q_{41}(l) & =q_{42}(l)=q_{43}(l)=q_{45}(l)=q_{46}(l)=q_{47}(l)=0 \\ 
q_{52}(l) & =q_{53}(l)=q_{54}(l)=q_{56}(l)=q_{57}(l)=q_{58}(l)=0 \\ 
q_{61}(l) & =q_{63}(l)=q_{64}(l)=q_{65}(l)=q_{67}(l)=q_{68}(l)=0 \\ 
q_{71}(l) & =q_{72}(l)=q_{74}(l)=q_{75}(l)=q_{76}(l)=q_{78}(l)=0 \\ 
q_{81}(l) & =q_{82}(l)=q_{83}(l)=q_{85}(l)=q_{86}(l)=q_{87}(l)=0%
\end{array}%
$
\end{example}

where $\mu =\frac{c}{b}$.

\begin{example}
Setting $n=6$, for $k=1,\;l=-3,\;a=i,\;b=i+1$ and $c=i-1$ $%
\begin{array}{lll}
J_{1} & = & diag(\lambda _{1},\lambda _{2},\lambda _{3},\lambda _{4},\lambda
_{5},\lambda _{6}) \\ 
& = & diag(-1.548i,-0.763i,0.371i,1.629i,2.763i,3.548i)%
\end{array}%
$
\end{example}

and%
\begin{equation*}
\begin{array}{ll}
(A_{1})^{-3} & =P_{1}(J_{1})^{-3}(P_{1})^{-1}=Q(-3)= \\ 
(q_{ij}(-3)) & =\left[ 
\begin{array}{ccc}
5i & -1.286-1.286i & -4.571 \\ 
1.286-1.286i & 0.429i & 1.571+1.571i \\ 
4.571 & -1.571+1.571i & 3.286i \\ 
-2.857-2.857i & 1.714 & 2.429-2.429i \\ 
-2.857i & 1.143+1.143i & 1.714 \\ 
-4+4i & -2.857i & -2.857-2.857i%
\end{array}%
\right. 
\end{array}%
\end{equation*}%
\begin{equation*}
\left. 
\begin{array}{cccc}
\qquad \qquad \qquad \qquad  & 2.857-2.857i & -2.857i & 4+4i \\ 
\qquad \qquad \qquad \qquad  & 1.714 & -1.143+1.143i & -2.857i \\ 
\qquad \qquad \qquad \qquad  & -2.429-2.429i & -1.714 & 2.857-2.857i \\ 
\qquad \qquad \qquad \qquad  & 3.286i & 1.571+1.571i & -4.571 \\ 
\qquad \qquad \qquad \qquad  & -1.571+1.571i & 0.429i & -1.286-1.286i \\ 
\qquad \qquad \qquad \qquad  & 4.571 & 1.286-1.286 & 5i%
\end{array}%
\right] .
\end{equation*}%
For $k=2,\;l=-4,\;a=3,\;b=1$ and $c=6$ 
\begin{equation*}
\begin{array}{lll}
J_{1} & = & diag(\lambda _{1},\lambda _{1},\lambda _{2},\lambda _{2},\lambda
_{3},\lambda _{3}) \\ 
& = & diag(-0.464,-0.464,3,3,6.464,6.464)%
\end{array}%
\end{equation*}%
and%
\begin{equation*}
\begin{array}{ll}
(A_{2})^{-4} & =P_{2}(J_{2})^{-4}(P_{2})^{-1}=Q(-4)= \\ 
(q_{ij}(-4)) & =\left[ 
\begin{array}{ccc}
5.395 & 0 & -3.112 \\ 
0 & 5.395 & 0 \\ 
-18.667 & 0 & 10.778 \\ 
0 & -18.667 & 0 \\ 
32.296 & 0 & -18.667 \\ 
0 & 32.296 & 0%
\end{array}%
\right. 
\end{array}%
\end{equation*}%
\begin{equation*}
\left. 
\begin{array}{cccc}
\qquad \qquad \qquad \qquad \qquad \qquad \qquad \qquad  & 0 & 0.897 & 0 \\ 
\qquad \qquad \qquad \qquad \qquad \qquad \qquad \qquad  & -3.112 & 0 & 0.897
\\ 
\qquad \qquad \qquad \qquad \qquad \qquad \qquad \qquad  & 0 & -3.112 & 0 \\ 
\qquad \qquad \qquad \qquad \qquad \qquad \qquad \qquad  & 10.778 & 0 & 
-3.112 \\ 
\qquad \qquad \qquad \qquad \qquad \qquad \qquad \qquad  & 0 & 5.395 & 0 \\ 
\qquad \qquad \qquad \qquad \qquad \qquad \qquad \qquad  & -18.667 & 0 & 
5.395%
\end{array}%
\right] .
\end{equation*}

\section{Complex Factorization}

The well-known Fibonacci polynomials $F(x)=\left\{F_{n}(x)\right\}
_{n=1}^{\infty }$ are defined in \cite{Thomas Koshy} by the recurrence
relation,

\begin{equation*}
F_{n}(x)=xF_{n-1}(x)+F_{n-2}(x) \;\; \left(n\geq 3\right)
\end{equation*}
where $F_{0}(x)=0$ , $F_{1}(x)=1$ and $F_{2}(x)=x$.

\begin{corollary}
Let the matrix $A_{k}$ be $n-$square ($n=2ks$, $s\in \mathbb{N}$) $k-$%
diagonal Toeplitz matrix as in (\ref{1}) with $a:=\mathit{x},b:=\mathit{i}$
and $c:=\mathit{i}$ \ where $\mathit{i}=\sqrt{-1}.$ Then 
\begin{equation}
\det (A_{k})=(F_{\frac{n}{k}+1}(x))^{k}  \label{26}
\end{equation}
where $F_{n}\left( x\right) $ are $n-$th Fibonacci polynomial.
\end{corollary}

\end{document}